\documentclass[12pt,twoside]{article}
\usepackage{latexsym, amssymb,amsthm,amsmath,lscape,epsfig,setspace,tikz,slashed,mathtools}
\usepackage[all]{xy}
\usepackage[retainorgcmds]{IEEEtrantools}
\usetikzlibrary{decorations.markings}
 \usetikzlibrary{arrows}
\usepackage[pdftex]{hyperref}
\usepackage[margin=10pt,font=small,labelfont=bf]{caption}
\usepackage{mathrsfs}

\usepackage{combelow}

\theoremstyle{plain}
\newtheorem{theorem}{Theorem}[section]
\newtheorem{proposition}[theorem]{Proposition}
\newtheorem{lemma}[theorem]{Lemma}
\newtheorem{corollary}[theorem]{Corollary}

\newtheoremstyle{namedthm}{}{}{\itshape}{}{\bfseries}{\!\!.}{0.5em}{\thmnote{#3 }}
\theoremstyle{namedthm}
\newtheorem*{namedtheorem}{Theorem}

\theoremstyle{definition}
\newtheorem{definition}[theorem]{Definition}

\newtheorem{example}[theorem]{Example}

\newenvironment{renumerate}%
{%
\begin{enumerate}}%
{\end{enumerate}%
}%

\newenvironment{remark}%
{\vskip6pt%
\noindent%
{\it Remark.}}%
{\vskip6pt}

{\vskip6pt%
\noindent%
{\it Remarks}. %
\begin{renumerate}}%
{\end{renumerate}\vskip6pt}

\makeatletter
\def\Ddots{\mathinner{\mkern1mu\raise\p@
\vbox{\kern7\p@\hbox{.}}\mkern2mu
\raise4\p@\hbox{.}\mkern2mu\raise7\p@\hbox{.}\mkern1mu}}
\makeatother

\newcommand{\h}{\mathfrak{h}}
\newcommand{\ka}{\kappa}
\newcommand{\B}{\mathbb{B}}
\newcommand{\diffto}{\xrightarrow{\raisebox{-0.2 em}[0pt][0pt]{\smash{\ensuremath{\sim}}}}}

\newcommand{\piff}{\pi_{\mathrm{aff}}}

\newcommand{\R}{\text{${\mathbb R}$}}

\renewcommand{\tilde}{\widetilde}

\newcommand{\e}{\text{$\varepsilon$}}

\newcommand{\supp}{\mathrm{supp}\,}

\newcommand{\mc}[1]{\text{$\mathcal{#1}$}}

\newcommand{\noqed}{\let\qed\relax}

\newcommand{\eps}{\varepsilon}

\date{} \usepackage{color} \definecolor{tocolor}{rgb}{.1,.1,.5}
\definecolor{urlcolor}{rgb}{.2,.2,.6}
\definecolor{linkcolor}{rgb}{.1,.1,.6}
\definecolor{citecolor}{rgb}{.6,.2,.1}
\hypersetup{backref=true, colorlinks=true, urlcolor=urlcolor,
  linkcolor=linkcolor, citecolor=citecolor}
\definecolor{darkgreen}{rgb}{0.0, 0.5, 0.0}

\begingroup
\hypersetup{linkcolor=blue}
\endgroup

\numberwithin{equation}{section}

\begin{document}

\title{Poisson structures with compact support}
\author{Gil R. Cavalcanti and 
 Ioan M\u{a}rcu\cb{t}
}
\maketitle

\abstract{We explicitly construct several Poisson structures with compact support. 
For example, we show that any Poisson structure on $\R^n$ with polynomial coefficients of degree at most two can be modified outside an open ball, such that it becomes compactly supported. We also show that a symplectic manifold with either contact or cosymplectic boundary admits a Poisson structure which vanishes to infinite order at the boundary and agrees with the original symplectic structure outside an arbitrarily small tubular neighbourhood of the boundary.
As a consequence, we prove that any even-dimensional manifold admits a Poisson structure which is symplectic outside a codimension one subset. 
}
\vskip12pt



\section{Introduction}

Already since Alan Weinstein's pioneering work \cite{MR723816}, a lot of the research in Poisson geometry has been focused on understanding ``the local structure of Poisson manifolds'', i.e., to prove local normal form theorems for Poisson structures around points or around certain submanifolds (see, e.g., \cite{MR3128977,MR4000576,MR794374,MR2776372,MR3005059,MR2178041,MR3632892,MR3250302,MR2144250,MR3261013}.)
These results deal with the question of how much of the local structure is determined by infinitesimal data. 
%
%
%
%

A natural follow-up problem is to understand how much the germ of a Poisson structure around a point or a submanifold influences its global behaviour. For example, in symplectic geometry, a germ of a symplectic structure might not extend to a global symplectic structure, because of various obstructions coming from symplectic topology. However, if we allow the extension to be a Poisson structure, no such restrictions are known. A concrete version of this problem, which we study in this paper in several settings is the following:\\[-6pt]

\begin{namedtheorem}[Poisson Extension Problem] \label{question}
Let $F$ be a closed subset of a manifold $M$, and $U$ be an open neighbourhood of $F$. Given a Poisson structure $\pi$ on a neighbourhood of $F$, does there exist a global Poisson structure $\tilde{\pi}$ on $M$, with the same germ as $\pi$ at $F$ and with $\mathrm{supp}(\tilde{\pi})\subset U$?
\end{namedtheorem}

This problem tests the flexibility of the Poisson relation, and plays an important role for developing \emph{h}-principles in Poisson geometry \cite{Pedro}.

We do not know of any situation where such an extension does not exist -- even if the condition on the support is dropped. To our knowledge, even in the case of when $F$ is just a point is an open question. When $F$ is a \emph{regular} point, a construction appeared in Waldmann's book \cite{Waldmann}. Inspired by his techniques, we obtain a positive answer to the Poisson Extension Problem in the case of a Poisson manifold with boundary for which the Poisson structure has a specific homogenous behaviour near the boundary  (see Theorem \ref{theorem:homogeneous}).  A simple consequence of this result generalises Waldmann's construction:

\begin{namedtheorem}[Lemma \ref{lemma:constant rank v2}] 
Let $\pi$ be a non-zero Poisson structure on $\R^{n}$ whose coefficients are polynomials of degree two. Then there exists a Poisson structure $\tilde{\pi}$ on $\R^n$ with $\mathrm{supp}(\tilde{\pi})=\B_1^n$ and $\pi|_{\B_{1/2}^n}=\tilde{\pi}|_{\B_{1/2}^n}$, where $\mathbb{B}^n_r$ denotes the closed ball of radius $r$ in $\R^n$. Moreover, if the coefficients of $\pi$ have degree one, then there exists a Poisson diffeomorphism $\phi: (\mathring{\mathbb{B}}_1^n,\tilde{\pi})\diffto (\R^n,\pi)$.
\end{namedtheorem}

The existence of such structures implies that the knowledge of the Poisson structure in a compact region may have no implications for the global Poisson geometry of the manifold. A direct consequence of Lemma \ref{lemma:constant rank v2} is that by choosing a  triangulation of a manifold $M$, one can endow the image of each simplex with a Poisson structure of arbitrary constant rank in the interior of that simplex and glue these together to produce a ``patchwork'' smooth Poisson structure on $M$ (see Corollary \ref{cor:patchwork}). In particular, every even dimensional manifold  admits a Poisson structure which is symplectic outside a codimension one subset. Here is a more immediate consequence of Lemma \ref{lemma:constant rank v2}:

\begin{namedtheorem}[Corollary \ref{cor:lie algebras}]
Given any $n$-dimensional Lie algebra $\mathfrak{g}$, there is a Poisson structure on $\R^n$, which is supported in $\B_1^n$ and whose restriction to $\mathring{\B}_1^n$ is globally isomorphic to the linear Poisson structure on $\mathfrak{g}^*$.
\end{namedtheorem}

We construct further examples by taking products:

\begin{namedtheorem}[Product Theorem]
For $i=1,2$, let $(M_i,\pi_i)$ be a manifold with a Poisson structure with compact support and let $U\subset M_1\times M_2$ be an open set containing $\supp(\pi_1)\times \supp(\pi_2)$. Then there is a Poisson structure $\Pi$ on $M_1 \times M_2$ such that
\begin{itemize}
\item 
$\Pi|_{\supp(\pi_1)\times \supp(\pi_2)} = \pi_1 + \pi_2$ and 
\item $\supp(\Pi)  \subset U$.
\end{itemize}
\end{namedtheorem}

Products allow us to move from balls to tubular neighbourhoods of manifolds:

{
\begin{namedtheorem}[Corollary \ref{cor:first extension}]
Let $(F,\pi_F)$ be a compact Poisson manifold and $\pi$ a Poisson structure on $\R^n$ with coefficients of degree at most two. The product $F \times \R^n$ admits a Poisson structure $\Pi$ with $\supp(\Pi)  \subset F \times \mathbb{B}_2^{n}$  such that 
\[\Pi|_{F\times \mathbb{B}_{1/2}^{n}}=\pi_F + \pi.\]
\end{namedtheorem}
}




{
As a more specific consequence, we obtain a positive answer to the  Poisson Extension Problem for symplectic germs around certain submanifolds: 

\begin{namedtheorem}[Corollary \ref{corollary:symplectic:submanifolds}]
Let $(M^{2n},\omega)$ be a symplectic manifold and $F\subset M $ be a compact symplectic submanifold whose normal bundle is symplectically trivial. For any neighbourhood $U$ of $F$, there is a Poisson structure $\tilde{\pi}$ on $M$ supported in $U$ which has the same germ around $F$ as $\omega$. 
\end{namedtheorem}
}

{
In the second part of the paper we consider the  Poisson Extension Problem for a compact domain endowed with a symplectic structure. More precisely, we take $(M, \omega)$ a compact symplectic manifold with boundary, and we try to extend $\pi=\omega^{-1}$ to a compactly supported Poisson structure on $M\cup_{\partial M} \partial M \times [0,1)$. The answer depends on the type of boundary $M$ has, and is the content of the {\it Poisson Extension Theorem}. The statement of this result needs the notions of {\it convexity} and {\it regular Pfaffian distributions}, so we will postpone to Section \ref{sec:extension}. Special cases of the theorem include boundary of contact type and of cosymplectic type. In particular, the former implies a positive answer to the Poisson Extension Theorem for symplectic germs around Lagrangian submanifolds. Another corollary is that every orientable compact four-manifold $M$ admits a Poisson structure which is symplectic outside a compact, codimension one submanifold, along which it vanishes to infinite order. 

}

These theorems are proved by explicit constructions. The main tool used in the proof of the Poisson Extension Theorem is the pure spinor description of Dirac and Poisson structures which  is common in generalized complex geometry but has not been exploited to its full potential within  the Poisson realm.

{ In Section \ref{section:homogeneous}, we prove the extension theorem for Poisson structures that have a certain homogeneity property at the boundary.}
We prove the existence of Poisson structures with support on the ball in Section \ref{sec:generalisations} and the Product Theorem in Section \ref{sec:products}. Sections \ref{sec:k-contact} and \ref{sec:convexity} introduce the concepts necessary to state and prove the Poisson Extension Theorem, which we do in Section \ref{sec:extension}. In Section \ref{sec:applications} we discuss the main applications of the results.
 
\subsection*{Acknowledgements} The authors would like to thank Marius Crainic and Pedro Frejlich for many discussions along the years, which have inspired and influenced the problems studied in this paper. 

In the first version of the paper, Lemma \ref{lemma:constant rank v2} was only for constant Poisson structures. Stefan Waldmann made us aware that this statement had appeared in his book \cite{Waldmann}. {His proof was much simpler than ours, and allowed us to extend the result to Poisson structures with polynomial coefficients of degree at most two.} We warmly thank him for his very useful feedback.

We also thank Alan Weinstein for his interest in the first version of the paper which led us to rephrase and re-structure some of the results.

Gil Cavalcanti thanks the Mathematics Department in Radboud University Nijmegen for their hospitality during the first stages of the project, and Ioan M\u{a}rcu\cb{t} thanks Instituto de Matem\'{a}tica Pura e Aplicada (IMPA) for its hospitality during the final stages of the project.

{
\section{Poisson structure homogeneous at the boundary}\label{section:homogeneous}

The following extension result turns out to be applicable in many situations.\\

\begin{theorem}\label{theorem:homogeneous}
Let $(M,\pi)$ be a compact Poisson manifold with boundary. Assume that there exist a vector field $X$ pointing outwards along $\partial M$, such that, around $\partial M$, $\pi$ decomposes as a finite sum of eigenvectors of $\mc{L}_X$:
\[\pi=\pi_1+\ldots+\pi_k,\quad \mc{L}_X\pi_i=\lambda_i\pi_i.\]
If $\lambda_i\leq 0$, then $\pi$ extends to a Poisson structure $\tilde{\pi}$ on $M\cup_{\partial M} \partial M \times [0,\infty)$ with compact support.
\end{theorem}

\begin{proof}
The vector field $X$ induces a collar neighbourhood $\mc{U}\simeq \partial M\times (-\epsilon,0]$ in which $X=\partial_t$. Then $\mc{L}_X\pi_i=\lambda_i\pi_i$ is equivalent to $\pi_i$ having the following form:
\[\pi_i=e^{\lambda_it}(\partial_t\wedge V_i+W_i),\]
with $V_i\in \mathfrak{X}^1(\partial M)$ and $W_i\in \mathfrak{X}^2(\partial M)$. We can use this formula to extend $\pi$ to a Poisson structure on $M\cup_{\partial M} \partial M \times [0,\infty)$. Let $f:[0,1)\to [0,\infty)$ be a diffeomorphism satisfying: 
	\begin{equation}\label{eq:f}
	f(t) = \begin{cases}
	t & \mbox { if }0\leq  t \leq \delta,\\
	e^{\frac{1}{1-t}}&  \mbox { if }  1-\delta \leq t< 1,
	\end{cases}
	\end{equation}
for some small $\delta>0$. Then $f$ induces a diffeomorphism 
\[\phi: M\cup_{\partial M}\partial M\times [0,1)\diffto M\cup_{\partial M}\partial M\times [0,\infty),\]
which is the identity on $M$, and is $\mathrm{Id}_{\partial M}\times f$ outside of $M$. We have that: 
\[\phi^*(\pi_i)=e^{\lambda_i f(t)}\Big(\frac{1}{f'(t)}\partial_t\wedge V_i+W_i\Big).\]
If $\lambda_i<0$ then this bivector extends smoothly, as zero to $\partial M\times [1,\infty)$. If $\lambda_i=0$, then the bivector extends smoothly as $W_i$ to $\partial M\times [1,\infty)$. We conclude that $\phi^*(\pi)$ extends smoothly to a Poisson structure, denoted $\pi_0$, which on $\partial M\times [1,\infty)$ it is given by $W=\sum_{\lambda_i=0}W_i$. 

Let $g:[0,\infty)\to [0,1]$ be a smooth function such that 
\[g(t)=1 \textrm{ for } t\leq 1,\quad  g(t)>0 \textrm{ for } t\in (1,2)\quad  \textrm{ and }\quad g(t)=0  \textrm{ for } t\geq 2.\] We regard $g$ as a function on $M\cup_{\partial M} \partial M\times [0,\infty)$, by setting $g=1$ on $M$. Note that $g$ is a Casimir function for $\pi_0$. So $\tilde{\pi}:=g\cdot \pi_0$ is a Poisson structure on $M \cup_{\partial M} \partial M\times [0,\infty)$, which is compactly supported and extends $\pi$. 
%
%
%
\end{proof}

As an immediate consequence of the theorem, we obtain:

\begin{corollary}\label{corollary:contact:and:cosymplectic}
Let $(M,\omega)$ be a compact symplectic manifold with boundary. If the boundary of $M$ is either of contact type or of cosymplectic type, i.e., if there exist a vector field $X$ pointing outwards along $\partial M$, such that, around $\partial M$, it satisfies: 
\[\textrm{either}\qquad  \mc{L}_X\omega=\omega \qquad  \textrm{or}\qquad \mc{L}_X\omega=0,\]
then $\omega$ extends to a Poisson structure $\tilde{\pi}$ on $M\cup_{\partial M} \partial M \times [0,\infty)$ with compact support.
\end{corollary}

\begin{remark}
The corollary will be generalized to a larger class of symplectic manifolds with boundary in Section \ref{sec:extension}. For this, let us note the two key steps of the proof of the theorem are:
\begin{enumerate}
\item Homogeneity allows us extend the Poisson structure $\pi$ on the entire $\partial M\times [0,\infty)$, with controlled growth;
\item After we pull back $\pi$ by a fast growing diffeomorphism, the Poisson structure extends to a constant Poisson structure.
\end{enumerate}
In the symplectic setting, the concept that will allow us to find an extension with `controlled growth' is that of {\it convexity}. Once this extension is produced, we can pull it back by a diffeomorphism as above and the concept that will ensure that the results structure can be further deformed by a Casimir is that of a {\it regular Pfaffian distribution}. Introducing these two concepts and stating the corresponding key steps will occupy Sections \ref{sec:k-contact} and \ref{sec:convexity}.
\end{remark}

}

{
\section{Poisson structures with support on the ball}\label{sec:generalisations}

%

The construction of Poisson structure supported on the ball follows easily from Theorem \ref{theorem:homogeneous}.

\begin{lemma}\label{lemma:constant rank v2}
Let $\pi$ be a non-zero Poisson structure on $\R^{n}$ whose coefficients are polynomials of degree two. Then there exists a Poisson structure $\tilde{\pi}$ on $\R^n$ with $\mathrm{supp}(\tilde{\pi})=\B_1^n$ and $\pi|_{\B_{1/2}^n}=\tilde{\pi}|_{\B_{1/2}^n}$, where $\mathbb{B}^n_r$ denotes the closed ball of radius $r$ in $\R^n$. Moreover, if the coefficients of $\pi$ have degree one, then there exists a Poisson diffeomorphism $\phi: (\mathring{\mathbb{B}}_1^n,\tilde{\pi})\diffto (\R^n,\pi)$.
\end{lemma}

\begin{proof}
Let $X=\sum_{i=1}^n x_i\partial_{x_i}$ denote the Euler vector field on $\R^n$. Then $X$ is transverse to the spheres $\partial B^n_{r}$. Decompose $\pi=\pi_0+\pi_1+\pi_2$, where $\pi_i$ has as coefficients of homogeneous polynomials of degree $i$. Then $\mc{L}_X \pi_i=(i-2)\pi_i$. We can therefore apply Theorem  \ref{theorem:homogeneous} to $M=\B_{1}^n$ and obtain the compactly supported extension $\tilde{\pi}$ of $\pi$. Let us look a bit closer to the proof. Note that, under the identification $X=\partial_t$, we have that $M\cup_{\partial M} \partial M\times [0,T)$ corresponds to $\B^n_{T+1}$. So if $\pi_2\neq 0$, $\mathrm{supp}(\tilde{\pi})=\B^n_{3}$. If $\pi_2=0$, then $\mathrm{supp}(\tilde{\pi})=\B^n_{2}$ and the construction gives also a Poisson diffeomorphism $\phi:(\mathring{\B}^n_2,\tilde{\pi})\diffto (\R^n,\pi)$. Finally, the radii can be easily fixed as in the statement.
\end{proof}
}

{
\begin{remark}
The construction from the affine case uses the same ideas as the proof of Proposition 4.1.20 in the book \cite{Waldmann}. More precisely, it uses the pullback of the Poisson structure to a ball via a rapidly growing diffeomorphism, as in Exercise 4.6 in \cite{Waldmann}. According to Waldmann \cite{Waldmann}, this construction was suggested by Weinstein, and indeed, for vector fields it is mentioned on page 66 in \cite{MR873455}. 
\end{remark}
}

\section{Products} \label{sec:products}

Now that we have established the existence of interesting Poisson structures with compact support, one can try and construct new examples from existing ones by taking products.

Of course, if $(M_1,\pi_1)$ and $(M_2,\pi_2)$ are manifolds endowed with Poisson structures with compact support, the product has also the product structure, $\pi_1 + \pi_2$ (where we have omitted the pullback maps), whose support is
\[(\mathrm{\supp}(\pi_1) \times M_2 )\cup (M_1 \times \mathrm{\supp}(\pi_2)).\]
If, say, $M_1$  is not compact and $\pi_2$ is nonzero, this support is not compact. As we will see next, this little problem can easily be fixed.

\begin{namedtheorem}[Product Theorem]
For $i=1,2$, let $(M_i,\pi_i)$ be a manifold with a Poisson structure with compact support and let $U\subset M_1\times M_2$ be an open set containing $\supp(\pi_1)\times \supp(\pi_2)$. Then there is a Poisson structure $\Pi$ on $M_1 \times M_2$ such that
\begin{itemize}
\item 
$\Pi|_{\supp(\pi_1)\times \supp(\pi_2)} = \pi_1 + \pi_2$ and 
\item $\supp(\Pi)  \subset U$.
\end{itemize}
\end{namedtheorem}

%
%
%
\begin{proof}
Since $\supp(\pi_i)$ is compact, there are open sets $U_i$ with $\supp(\pi_i) \subset U_i$ and $U_1 \times U_2 \subset U$. By the existence of partitions of unity, there are smooth functions $\chi_i\colon M_i \to [0,1]$ such that $\chi_i|_{\supp(\pi_i)} = 1$ and $\supp(\chi_i) \subset U_i$. Notice that $\chi_i$ is a Casimir function for $\pi_i$. 
We endow  $M_1 \times M_2$ with the bivector field 
\[\Pi = \chi_2 \pi_1+\chi_1\pi_2,\]
which is clearly a Poisson structure and satisfies the conditions from the theorem.
\end{proof}

\section{The geometry of hyperplane distributions}\label{sec:k-contact}


Let $\kappa$ be a hyperplane distribution on a manifold $N$. Locally one can always find a 1-form, $\gamma$,  whose kernel is $\kappa$ and if $\kappa$ is coorientable there is a global 1-form defining it. The degree of integrability of $\kappa$ is governed by the restriction of $d\gamma$ to $\kappa$. The case when the restriction has constant rank deserves special attention.

\begin{definition}[See \cite{MR1083148,MR882548}]
 A \emph{Pfaffian distribution of class} $2k+1$ is a hyperplane distribution $\kappa$ such that  $\gamma \wedge(d\gamma)^{k} \neq  0$ and $\gamma \wedge(d\gamma)^{k+1}= 0$, where $\gamma$ is any local 1-form defining $\kappa$. 
Further, $\kappa$ is  a {\it regular Pfaffian distribution (of class $2k+1$)} if $\gamma \wedge(d\gamma)^k$ is nowhere zero.
\end{definition}

Note that, if we replace in $\gamma\wedge (d \gamma)^p$ the form $\gamma$ by a multiple $\varphi\gamma$, we obtain $\varphi^{p+1}\gamma\wedge (d\gamma)^p$. Therefore, the definition can be checked using any 1-form $\gamma$ defining $\kappa$.

\begin{example}
\begin{itemize}
\item
A contact structure $\ka$ on $N^{2k+1}$ is the same as a regular Pfaffian distribution of maximal class $2k+1$. 
\item At the other extreme, if $\kappa$ is a regular Pfaffian distribution of class 1, then it is in fact integrable as $\gamma \wedge d\gamma \equiv 0$ is precisely the integrability condition for $\kappa$. Therefore, the larger the values of $k$, the more nonintegrable $\ka$ is.

\item Other intermediate examples can be obtained as follows. Let $\phi:N\to B$ be a surjective submersion, and let $\ka_B$ be a contact structure on $B$. Then $\ka:=(d\phi)^{-1}(\ka_B)$ is a regular Pfaffian distribution of class $2k+1=\dim(B)$.
\end{itemize}
\end{example}

A classical theorem of Darboux says that this last example provides in fact the local structure of all regular Pfaffian distributions: 

\begin{theorem}[Darboux \cite{Darboux1}]
Let $\kappa$ be a regular Pfaffian distribution of class $2k+1$ on $N^n$. Then $N$ can be covered by charts with coordinates $(x_1,\ldots,x_{n})$, in which 
\[\ka=\ker \big(dx_1+x_2d x_3+\ldots + x_{2k}d x_{2k+1}\big).\]
\end{theorem}

This theorem has also the following global interpretation:

\begin{proposition}
Let $\kappa$ be a regular Pfaffian distribution of class $2k+1$ on $N^n$. Then 
\begin{equation}\label{eq:integrable distribution}
\h:= \ker (d \gamma|_{\ka})\subset TN
\end{equation}
is an involutive distribution of rank $n-2k-1$, where $\gamma$ is any local 1-form defining $\kappa$. Moreover, $\ka$ induces a transverse contact structure to the foliation generated by $\h$, i.e.
\begin{enumerate}
\item[1.] Any transversal $X\subset N$ to $\h$ has an induced contact structure $\ka_X:=TX\cap \ka$;
\item[2.] Any holonomy transformation is a contactomorphism 
$\mathrm{hol}:(X_1,\ka_{X_1})\diffto (X_2,\ka_{X_2})$.
\end{enumerate}
In particular, if $\h$ gives rise to a simple foliation, i.e., if there is a surjective submersion with connected fibers, $\phi:N\to B$, with $\h=\ker d \phi$, then the base $B$ carries a contact structure $\ka_B$ such that $\ka=(d\phi)^{-1}(\ka_B)$. 
\end{proposition}

The involutive distribution $\h$ will be called the \emph{kernel of the regular Pfaffian structure}.

\begin{proof} All these statements are by proved picking appropriate Darboux charts. First, involutivity of $\h$ follows because, in Darboux coordinates, $\h$ is spanned by $\partial_{x_{2k+2}},\ldots, \partial_{x_{n}}$. 

Similarly, for {\it 1.}, the first coordinates, $(x_1,\ldots,x_{2k+1})$ form a local coordinate system on $X$, in which $\ka|_X$ becomes the standard contact structure on $\R^{2k+1}$. 

Item {\it 2.} follows because $\h$ normalizes $\ka$:
\[X\in \Gamma(\h),\quad Y\in \Gamma(\ka)\implies [X,Y]\in \Gamma(\ka)\]
as can be easily checked using a local Darboux chart. 

The last claim follows from {\it 1.} and {\it 2.}
\end{proof}
%
%
%

The property of regular Pfaffian distributions which is relevant to us is the following interpolation construction. For this, fix a smooth function $f\colon [0,1) \to \R$ such that
\begin{equation}\label{eq:f}
f(t)=e^{e^{\frac{1}{1-t}}}, \quad \mbox{ for all } t \mbox{ close to } 1.
\end{equation}
The key property of $f$ is that, for $i\geq 1$, the functions $\frac{1}{f^i}$, $\frac{1}{f' f^i}$ and $\frac{f}{f'}$ all have zeros of infinite order at $t=1$, i.e., their extension by $0$ for $t\geq 1$ is smooth.


\begin{namedtheorem}[Dirac Interpolation Lemma]
Let $N$ be a manifold with a regular Pfaffian distribution $\kappa$, defined by a 1-form $\gamma$ and with kernel $\h$. Then the Dirac structure $D$ on $N \times [0,1)$ with spinor line spanned by
\begin{equation}\label{eq:Dirac}
\rho = e^{d(f(t) \gamma)},
\end{equation}
extends smoothly to $N\times [0,1]$. Further, $D$ agrees with the Dirac structure $ \h \oplus \h^0$ to infinite order at $t=1$, where $\h^0 \subset T^*(N\times [0,1])$ is the annihilator of $\h$. Therefore $N \times \R $ admits a Dirac structure which is gauge equivalent to $T(N\times \R)$ for $t < 1$ and is $\h \oplus \h^0$ for $t \geq 1$.
\end{namedtheorem}
\begin{proof}
We need to study the behaviour of the spinor line of $D$ as $t$ approaches $1$. Notice that the particular spinor $\rho$ which defines $D$ on $N \times [0,1)$ blows up at $1$, hence it cannot be used to trivialise the spinor line beyond that point. So, to check if the line can be extended beyond that point we must renormalize it. First we compute $\rho$:
\begin{align*}
\rho &= e^{d(f(t) \gamma)} =e^{f' dt\wedge  \gamma + f d\gamma}\\
&=(1+ f' dt \wedge \gamma)\wedge\big(\sum_{j=0}^{k+1}\frac{1}{j!}f^j(d\gamma)^j\big)\\
& =1 + \sum_{j=0}^{k}f^j\frac{1}{j!}\Big(\frac{1}{j+1}f d\gamma + \frac{f'}{j!}dt\wedge \gamma\Big)\wedge (d\gamma)^{j} ,
\end{align*}
where we used that $\gamma \wedge (d\gamma)^{k+1} \equiv 0$ to conclude that $(d\gamma)^{k+2} \equiv 0$ and hence there are no further terms in the exponential expansion. We normalize $\rho$ by dividing it by $f'f^{k}$ to get
\[
\frac{\rho}{f' f^k} = \frac{1}{f' f^k} + \sum_{j=0}^{k} f^{j-k} \frac{1}{j!}\Big(\frac{1}{(j+1)}\frac{f}{f'}d\gamma + dt\wedge \gamma \Big)\wedge (d\gamma)^{j}.
\]
By our choice of $f$, $\frac{\rho}{f' f^k}$ extends smoothly at $t=1$, and equals $\frac{1}{k!} dt\wedge \gamma\wedge (d\gamma)^k$ to infinite order at $t=1$, which is the spinor line for $\h \oplus \h^0$.
\end{proof}

\begin{remark}
If $\kappa$ is a contact structure, then $\h=0$ and, by taking $f=0$ around $t=0$, we obtain that the Dirac structure $D$ interpolates between $T(N \times \R)$ and $T^*(N\times \R)$.
\end{remark}

\section{Convexity}\label{sec:convexity}

Here we use standard definitions of convexity of symplectic manifolds as presented, for example, in the survey \cite{MR1634561}.

The notion of convexity belongs to the world of almost complex manifolds. Given an almost complex manifold with boundary, $(M,J)$, we can consider the hyperplane distribution on $\partial M$ given by $\kappa = T \partial M \cap J (T\partial M)$. Alternatively, $\kappa$ can be described as the maximal $J$-invariant subspace of $T\partial M$.

\begin{definition}
The {\it Levi form} on $\partial M$ is the bilinear form
$$L\colon \kappa \times \kappa \to T(\partial M)/\kappa,\qquad L(V,W)=\mathrm{pr}( [JV,W]), \qquad \mbox{ for } V,W \in \Gamma(\kappa).$$
where $\mathrm{pr}\colon T(\partial M) \to T(\partial M)/\kappa$ is the natural projection.
\end{definition}
Notice that $L$ is completely determined by $J|_\kappa$, so any two almost complex structures that agree on $\kappa$ give rise to the same Levi form.

We coorient $\kappa$ by declaring that a vector $\mathrm{pr}(V) \in T(\partial M)/\kappa$ is positive if $JV$ points inwards. With this convention, we can talk about the sign of  $L(V,W)$.
\begin{definition}
An almost complex manifold $(M,J)$ has {\it convex} boundary if the Levi form is positive semi-definite, i.e., if
\[L(V,V)=\mathrm{pr}( [JV,V])\geq 0, \qquad \mbox{ for all } V\in \kappa.\]
\end{definition}

%

Since $\kappa$ is coorientable, there is a 1-form $\gamma \in \Omega^1(\partial M)$ such that $\ker(\gamma) = \kappa$ and $\gamma(V)> 0$ if and only if $\mathrm{pr}(V)> 0$. We call any such form a {\it positive primitive for the Levi form}. It follows that $L$ is positive semi-definite if and only if $\gamma([JV,W])$ is positive semi-definite and Cartan's formula relates the latter with the exterior derivative of $\gamma$:
$$ \gamma([JV,W]) =d\gamma(W,JV) \qquad \mbox{ for all } V,W \in \Gamma(\kappa).$$
So $\partial M$ is convex if and only if  
\[d\gamma(V,JV)\geq 0, \qquad \mbox{ for all } V\in \kappa.\]

The notion of convexity goes through to the symplectic world by choosing an almost complex structure which tames $\omega$:

\begin{definition}
Let $(M,\omega)$ be a symplectic manifold. The boundary of $M$  is {\it pseudoconvex} if there is an almost complex structure $J$ defined in a neighbourhood $U$ of $\partial M$ which \emph{tames} $\omega$, i.e., 
\[\omega(V,JV)>0, \qquad \mbox{ for all } V\in TU,\ V\neq 0\]
and with respect to which $\partial M$ is convex.
\end{definition}

We denote the restriction of $\omega$ to the boundary by
\[\sigma:=\omega|_{\partial M}\in \Omega^2(\partial M).\]
Notice that if $J$ tames $\omega$, then $\kappa = T\partial M \cap J (T\partial M)$ is transverse to the natural rank-one distribution on $\partial M$ given by $\mathrm{ker}\,\sigma $. Therefore, if we  let $\gamma \in \Omega^1(\partial M)$ be a positive primitive for the Levi form of $J$, by the Coisotropic Neighbourhood Theorem \cite{MR633290}, $\omega$, in a neighbourhood of $\partial M$, is equivalent to the following form 
defined on $\partial M\times (-\e,0]$: 
\begin{equation}\label{eq:omega_N}
\omega_{0} := \sigma + d(t \gamma).
\end{equation}

The relevance of convexity comes from the next lemma:

\begin{namedtheorem}[Symplectic Extension Lemma]
Let $(M,\omega)$ be a symplectic manifold with pseudoconvex boundary. Then the manifold $M \cup_{\partial M} (\partial M \times [0,\infty))$ has a symplectic structure which equals $\omega$ on $M$ and $\omega_0$ from \eqref{eq:omega_N} on $\partial M \times [0,\infty)$.
\end{namedtheorem}
\begin{proof}
We work on a tubular neighbourhood $U \simeq  {\partial M}\times (-\e,0]$ of $\partial M$ where the symplectic structure is given by \eqref{eq:omega_N}. We need to show that $\omega_0$ is non-degenerate on $\partial M\times [0,\infty)$. In the decomposition $T^*(\partial M)=\kappa^*\oplus \langle \gamma\rangle $, we have that $\sigma \in \Gamma(\wedge^2 \kappa^*)$ and we can write
\[
d\gamma = \alpha + \beta\wedge \gamma,
\]
with $\alpha \in \Gamma(\wedge^2\kappa^*)$ and $\beta \in  \Gamma(\kappa^*)$. Then the top power of $\omega_0$ is $n(\sigma+ t\alpha)^{n-1}\wedge \gamma \wedge dt$. So it suffices to show that, for $t\geq 0$, $\sigma+t\alpha$ is non-degenerate on $\ka$. This follows because $J|_{\ka}$ tames $\sigma|_{\partial M}$ and $\partial M$ is convex for $J$:
\[(\sigma+t\alpha)(V,JV)>0, \qquad \mbox{ for all } V\in \ka,\ V\neq 0.\qedhere\]

\begin{remark}
To keep the arguments short, we assumed that the whole boundary of $M$ is convex, but it is clear that one can apply the Symplectic Extension Lemma to individual components of the boundary as long as those are pseudoconvex. The same comment applies to the results of the next sections.
\end{remark}
\end{proof}

\section{Poisson Extension Theorem}\label{sec:extension}

Now we set up, state and prove our main extension result for symplectic structures.  We start by identifying the types of symplectic structures on manifolds with boundary that can be extended to Poisson structures with compact support.

\begin{definition}\label{def:conver k-contact}
We say that a symplectic manifold, $(M,\omega)$, has {\it pseudoconvex boundary of regular Pfaffian type} if there exists an almost complex structure $J$ such that 
\begin{itemize}
\item[1)] $\partial M$ is convex for $J$,
\item[2)]$\ka:=T\partial M\cap J(T\partial M)$ is a regular Pfaffian distribution, and
\item[3)] the restriction of $\omega$ to the kernel of $\kappa$ is non-degenerate.
\end{itemize}
\end{definition}

Because the space of almost complex structure taming $\omega$ is contractible, it suffices to assume the existence of $J$ around (or just along) $\partial M$.

As we will see in Section \ref{sec:applications}, symplectic manifolds with pseudoconvex boundary of regular Pfaffian type include symplectic manifolds with contact or with cosymplectic boundary. 

The definition can be mildly rephrased as follows.

\begin{lemma}\label{lem:alternative}
The boundary of a symplectic manifold, $(M^{2n},\omega)$, is pseudoconvex of regular Pfaffian type of class $2k+1$ if and only if there is an outwards pointing vector field, $X$, for which
\begin{enumerate}
\item[a)]  $\gamma: = i_X \omega|_{\partial M}$ defines a regular Pfaffian distribution of class $2k+1$, $\kappa$, on $\partial M$,
\item[b)] $\omega^{n-k-1}|_{\partial M}\wedge (d\gamma)^k \wedge \gamma$ is a volume form on $\partial M$, and
\item[c)] there is an almost complex structure $J_0$ on $\kappa$ such that:
\[\omega(V,J_0V)>0 \quad \textrm{and}\quad  (\mathscr{L}_X\omega)(V,J_0V)\geq 0,\quad \mbox{ for all } V\in \ka\backslash\{0\}.\]
\end{enumerate}
\end{lemma}
\begin{proof}
Assume that $\partial M$ is pseudoconvex of regular Pfaffian type. Let $X$ be a vector field pointing outwards. We may assume that, along $\partial M$, $X\in \ka^{\perp_{\omega}}$. This follows because $J$ tames $\omega$, and so $\ka$ is a symplectic subspace, and so we can change $X$ by a section of $\ka$ so that $X\in \ka^{\perp_{\omega}}$. Then for $\gamma:=i_X\omega|_{\partial M}$, we have $\ker\gamma=\ka$. To show that $\gamma$ is a positive primitive of the Levi form, we use the direct sum decomposition: 
\begin{equation}\label{eq:some:decomp}
TM|_{\partial M}=\ka\oplus \ker (\omega|_{\partial M})\oplus \langle X \rangle.
\end{equation}
Let $V$ be the section of $\ker (\omega|_{\partial M})$ such $\gamma(V)=1$. We can decompose $JV=A+aV+bX$, with $A\in \ka$ and $a,b\in \R$. We have that:
\[0<\omega(V,JV)=\omega(V,A+aV+bX)=b\omega(V,X)=-b\gamma(V)=-b.\]
Thus, $b<0$, which shows that $JV$ points inwards. 

Item a) is the same as item 2), and item b) is just an algebraic reformulation of item 3). For item c), we take $J_0:=J|_{\ka}$. Since $J$ tames $\omega$ and $(\mathscr{L}_X\omega)|_{\partial M}=d\gamma$, item c) follows from item 1). 

Conversely, assume that $X$ is given with properties a), b) and c). Item c) implies that $\omega|_{\ka}$ is symplectic, and therefore the decomposition \eqref{eq:some:decomp} holds. Let $V$ be the section of $\ker (\omega|_{\partial M})$ satisfying $\gamma(V)=1$. We define the almost complex structure $J$ by setting on $TM|_{\partial M}$
\[J|_{\ka}=J_0, \quad JV=-X, \quad J X=V,\]
and extending to $M$ so that it tames $\omega$. Clearly $T(\partial M)\cap J (T\partial M)=\ka$. Since $JV=-X$ points inwards, $\gamma$ is a positive primitive of the Levi of $J$. 

Arguments as above, immediatly show that: 2) $\Rightarrow$ a), 3) $\Rightarrow$ b), and 1) $\Rightarrow$ c). 
\end{proof}

\begin{remark}\label{re:on:boundary} Let $(M,\omega)$ be a symplectic manifold whose boundary is pseudoconvex of regular Pfaffian type. Then $\partial M$ carries a Poisson structure $\pi_{\partial M}$ with leaves given by the kernel of the regular Pfaffian distribution $\ka$ and the symplectic structure on the leaves is given by the restriction of $\omega$. Equivalently, the spinor line corresponding to $\pi_{\partial M}$ is spanned by:
\[e^{\sigma}\wedge \gamma\wedge (d\gamma)^k,\quad \textrm{where}\quad \sigma:=\omega|_{\partial M}.\]
\end{remark}

{
\begin{namedtheorem}[Poisson Extension Theorem]\label{theo:extension theorem}
Let $(M,\omega)$ be a compact symplectic manifold whose boundary is pseudoconvex of regular Pfaffian type. Then $\omega$ can be extended to a compactly supported Poisson structure on $M\cup_{\partial M}(\partial M\times [0,\infty))$, denoted $\tilde{\pi}$, which satisfies: 
\begin{itemize}
\item $\tilde{\pi}$ is symplectic on $M\cup_{\partial M}(\partial M\times [0,1))$;
\item $\tilde{\pi}$ is regular on $\partial M\times [1,2)$, and its symplectic leaves integrate the distribution $H\simeq \h\times 0_{[1,2)}$, where $\h$ is the kernel of the regular Pfaffian distribution;
\item $\tilde{\pi}$ vanishes on $\partial M\times [2,\infty)$.
\end{itemize}
\end{namedtheorem}
}


\begin{proof}
Let $J$ be the almost complex structure for which $\partial M$ is pseudo-convex with regular Pfaffian structure $\ka$, and let $\gamma$ be a positive primitive of the Levi form. By the Symplectic Extension Lemma, we find an arbitrary small collar $U_{\e}\simeq \partial M\times (-\e,0]$ around the $\partial M$, such that $\omega$ can be extended to a symplectic structure $\tilde{\omega}$ on $M \cup_{\partial M} (\partial M \times [0,\infty))$ satisfying: 
\[\tilde\omega|_{\partial M\times (-\e,\infty) }=\sigma+ d(t\gamma),\quad \textrm{where} \quad \sigma=\omega|_{\partial M}.\]
Let $f:[0,1)\diffto [0,\infty)$ be a diffeomorphism which is the identity near $t=0$ and is given by \eqref{eq:f} near $t=1$, and define the diffeomorphism
\[\phi:M \cup_{\partial M} (\partial M \times [0,1))\diffto M \cup_{\partial M} (\partial M \times [0,\infty)),\] \[\phi|_{M}=\mathrm{Id}_M,\quad \phi|_{\partial M\times [0,1)}=\mathrm{Id}_{\partial M}\times f.\]
The pullback $\phi^*(\tilde{\omega})$ is an extension of $\omega$ to a symplectic structure $M \cup_{\partial M} (\partial M \times [0,1))$. We analyse this form as a Dirac structure. On $\partial M \times [0,1)$, it is represented by the spinor:
\[e^{\phi^*\tilde{\omega}} = e^{\sigma + d(f(t)\gamma)} = e^\sigma e^{d(f(t)\gamma)}.
\]
The Dirac Interpolation Lemma implies that the 2-form $\phi^*(\tilde{\omega})$ extends as a Dirac structure to $M\cup_{\partial M} \partial M\times [0,\infty)$ with defining spinor given on ${\partial M}\times [1,\infty)$ by:
\[e^\sigma \wedge dt\wedge \gamma \wedge (d\gamma)^k.\]
By assumption, the top element, $\sigma^{n-k-1}\wedge dt \wedge \gamma \wedge (d\gamma)^k$, is a volume form, and so this Dirac structure is in fact Poisson $\pi_{0}$ on $M\cup_{\partial M} \partial M\times [0,\infty)$. Note that $\pi_0$ extends $\omega$, it is symplectic on $\partial M\times (0,1)$, and on $\partial M\times [1,\infty)$ it equals the Poisson structure $\pi_{\partial M}$ described in Remark \ref{re:on:boundary}. To obtain $\tilde{\pi}$, we multiply $\pi_0$ with a Casimir function, as in the proof of Theorem \ref{theorem:homogeneous}. 
%
%
\end{proof}

\section{Applications}\label{sec:applications}

\subsection{Poisson structures with small support}

A direct consequence of Lemma \ref{lemma:constant rank v2} is Proposition 4.1.20 \cite{Waldmann} (whose proof we generalized):
\begin{corollary}
For any $0\leq 2r\leq n$, there is a Poisson structure $\pi_{2r,n}$ on $\R^{n}$, which vanishes outside $\B_1^n$ is has constant rank equal to $2r$ on $\mathring{\B}_1^n$.
\end{corollary}

Recall that linear Poisson structures are in one-to-one correspondence with Lie algebras (see e.g., \cite{MR4328925}). So we obtain another obvious consequence of Lemma \ref{lemma:constant rank v2}, worth mentioning:
\begin{corollary}\label{cor:lie algebras}
Let $(\mathfrak{g},[\cdot,\cdot])$ be a Lie algebra and let $(\mathfrak{g}^*,\pi_{\mathfrak{g}})$ be the corresponding linear Poisson structure on its dual. There is a Poisson structure $\tilde{\pi}_{\mathfrak{g}}$ on $\mathfrak{g}^*$ which coincides with $\pi_{\mathfrak{g}}$ on the ball $\B_{1/2}$ and vanishes outside $\B_1$.
\end{corollary}

Another particular, interesting case is that of a product of a constant, non-degenerate Poisson structure with a linear Poisson structure. By using Lemma \ref{lemma:constant rank v2} for these, we obtain the following variation of Proposition 4.1.20 \cite{Waldmann}:

\begin{corollary}
Let $(M,\pi)$ be a Poisson manifold. For any $p\in M$, and any neighborhood $U$ of $p$ in $M$, there exists a Poisson structure $\tilde{\pi}$ with support in $U$, which has the same first jet as $\pi$ at $p$.
\end{corollary}

\begin{proof}
By Weinstein's Splitting Theorem \cite{MR723816}, there is a chart centred at $p$, in which $\pi$ is the product of the standard symplectic structure and a Poisson structure which vanishes at $0$. Then, the first order expansion of $\pi$ in that chart at $0$ is an affine Poisson structure $\piff$. By applying Lemma \ref{lemma:constant rank v2} to $\piff$, we obtain a Poisson structure $\tilde{\pi}$, with the same germ as $\piff$ at $p$, and which we may assume to be supported inside $U$.
\end{proof}

We remark that a general affine Poisson structure is determined by a Lie algebra endowed with a 2-cocycle (see e.g., Proposition 2.29 \cite{MR4328925}).

\subsection{Patchwork Poisson structures}

The ``bump Poisson structures'' of Lemma \ref{lemma:constant rank v2} vanish flatly at the boundary, and therefore they can be used to decorate any manifold with a patchwork of regular Poisson structures.

\begin{corollary}\label{cor:patchwork}
Given a smooth $n$-dimensional manifold,  $M$, with a smooth triangulation
\[\Sigma = \{\sigma_i\colon \triangle^n \to M\colon i \in I\},\]
and a collection of numbers $\{r_i\}_{i\in I}$ with $r_i \in  \{0,\dots, \lfloor \tfrac{N}{2}\rfloor\} $, there is a smooth Poisson structure on $M$ which  has (constant) rank $2r_i$ in the interior of $\sigma_i(\triangle^n)$ and vanishes at  lower dimensional simplices. 
\end{corollary}
\begin{proof}
It suffices to prove that for every $r \leq\lfloor \tfrac{n}{2}\rfloor$ there is a Poisson structure on an $n$-simplex, $\triangle^n$, which has rank $2r$ in $\mathring\triangle^n$ and vanishes to infinite order at the boundary, as the desired structure is obtained by patching a collection of such structures together.
  
To prove that an $n$-simplex admits such a Poisson structure we observe that for all $n$ there is a map $\phi_n \colon \triangle^n \to \mathbb{B}^n_1$ which is a diffeomorphism between  the interior of an $n$-dimensional simplex and the interior of the $n$-dimensional unit ball and whose singularities at the boundary involve only roots and polynomials.

One such map can be constructed inductively on the dimensions involved. For $n=1$, the interval $[-1,1]$ is both a 1-simplex and the unit ball and we can take the identity map as a diffeomorphism. For the inductive step, assume we have a map  $\phi_n \colon \triangle^n \to \mathbb{B}_1^n$ with the stated properties, where $\triangle^n \subset \R^n$ is a symplex and consider $\triangle^{n+1} \subset \R^{n}\times \R$ the pyramid with base the symplex $\triangle^n\times \{-1\} \subset \R^{n}\times \R$ and vertex at $(0,1) \in \R^n\times \R$:
\[\triangle^{n+1} = \left\{(\tfrac{1-t}{2}x,t)\in \R^{n} \times \R \colon x\in \triangle^n, t\in [-1,1]\right\}.\]  
Then define
\[ \phi_{n+1} \colon \triangle^{n+1} \to \mathbb{B}^{n+1},\qquad \phi_{n+1}(\tfrac{1-t}{2}x,t) = (\sqrt{1-t^2} \phi_N(x), t).\]
One can readily check that this map has the stated properties. For example, starting with $\triangle^1 = [-1,1]$ and $\phi_1 = \mathrm{Id}$ we have
\[\triangle^2 = \left\{(x,y) \in \R^2\colon  |y| \leq1, |x| \leq \frac{1-y}{2}\right\},\qquad   \phi_{2}(x,y) = \left(2\sqrt{\frac{1+y}{1-y}}\, x,y\right).\]

It follows that the Poisson structure  of rank $2r$ in $\mathring{\mathbb{B}}_1^{n}$ constructed in Lemma \ref{lemma:constant rank v2} can be pulled back to a Poisson structure on the simplex which has rank $2r$ in $\mathring{\triangle}^{n}$ and vanishes flatly at the boundary.
\end{proof}

\begin{corollary}
Every even dimensional smooth manifold admits a Poisson structure which is symplectic outside a codimension-one subset.
\end{corollary}

\begin{remark} The simplices in Corollary \ref{cor:patchwork} could have been decorated with any affine Poisson structures. For example, we could have chosen $N$-dimensional Lie algebras, $\{\mathfrak{g}_i\}_{i\in I}$, such that the resulting Poisson structure restricted to the simplex $\sigma_i(\mathring{\triangle}^N)$ is Poisson diffeomorphic to the linear Poisson structure $(\mathfrak{g}_i^*,\pi_{\mathfrak{g}_i})$, for all $i\in I$.
\end{remark}

\subsection{Submanifolds with trivial normal bundle}\label{ex:submanifolds}

The Product Theorem yields interesting examples already when one of the manifolds involved is compact.
\begin{corollary}\label{cor:first extension}
{
Let $(F,\pi_F)$ be a compact Poisson manifold and $\pi$ a Poisson structure on $\R^n$ with coefficients of degree at most two.}
Then $F \times \R^n$ admits a Poisson structure $\Pi$ such that 
\begin{itemize}
\item 
$\Pi|_{F \times \mathbb{B}_{1/2}^N}=\pi_F + \pi$,
\item $\supp(\Pi)  =(F \times \mathbb{B}_1^{N})\cup (\supp(\pi_F)\times \mathbb{B}_2^N)$.
\end{itemize}
\end{corollary}

{
We can further specialize the statement to symplectic structures, and obtain that the Poisson Extension Problem has a positive answer around certain symplectic submanifolds:
}

\begin{corollary}\label{corollary:symplectic:submanifolds}
Let $(M^{2n},\omega)$ be a symplectic manifold without boundary and $F\subset M $ be a compact symplectic submanifold whose normal bundle is symplectically trivial. For any neighbourhood $U$ of $F$, there is a Poisson structure $\pi$ on $M$ which has the same germ around $F$ as $\omega$ and is supported in $U$.
\end{corollary}

{

\begin{proof}
By the normal form around symplectic manifolds, a neighbourhood of $F$ is isomorphic to a product $F\times \B_{\e}^{2k}$, with the product symplectic structure. Applying Corollary \ref{cor:first extension}, the conclusion follows. 
Alternatively, we can apply Theorem \ref{theorem:homogeneous} for $M=F\times \B_{\e}^{2k}$ and the Euler vector field on $\R^{2k}$. Alternatively, the result also follows form the Poisson Extension Theorem: the boundary of $M$ is pseudoconvex of regular Pfaffian type of class $2k-1$. 
\end{proof}
}
{
\begin{remark}
By using the normal form theorem around Poisson transversals \cite{MR3632892}, similar arguments implies a positive answer to the Poisson Extension Problem around Poisson transversals with symplectically trivial normal bundle. 
\end{remark}

}

\subsection{Contact boundary}\label{ex:contact}

Let $(M^{2n},\omega)$ have boundary of contact type, that is, there is an outwards pointing vector field $X$ defined around $\partial M$ such that $\mc{L}_X \omega = \omega$. Then the boundary of $(M^{2n},\omega)$ is pseudoconvex of regular Pfaffian type of class $2n-1$. Indeed, $\gamma = i_X  \omega|_{\partial M}$ defines a contact structure on $\partial M$, i.e., a regular Pfaffian distribution of class $2n-1$ with kernel $\h=0$. Thus the first two conditions in Lemma \ref{lem:alternative} are satisfied. Further, since $\mc{L}_X \omega = \omega$, the third condition holds for any almost complex structure on $\ka:=\ker \gamma$  taming $\omega|_{\ka}$. 

{
Therefore, the Poisson Extension Theorem gives an alternative approach to the first part of Corollary \ref{corollary:contact:and:cosymplectic}:

\begin{corollary}\label{cor:contact}
Let $(M^{2n},\omega)$ be a symplectic manifold with contact boundary. Then $\omega$ extends to a Poisson structure $\tilde{\pi}$ on $M\cup_{\partial M} \partial M\times [0,\infty)$, which is symplectic on $M\cup_{\partial M} \partial M\times [0,1)$ and vanishes on the rest $\partial M\times [1,\infty)$.
\end{corollary}
}

For the unit ball $M:=\B^{2n}_1\subset \R^{2n}$ with the standard symplectic structure, the corollary recovers the constant, symplectic case of Lemma \ref{lemma:constant rank v2}. 

Another class of examples are obtained as follows. Let $L$ be a compact Riemannian manifold. Then the disk bundle in the cotangent bundle 
\[M_{\eps}:=\{\xi\in T^*L\,  : \, |\xi|\leq \e\},\quad \omega:=\omega_{\mathrm{can}}|_{M_{\e}},\]
is symplectic with contact boundary. Hence, by Corollary \ref{cor:contact}, $T^*L$ admits a compactly supported Poisson structure $\tilde{\pi}$ which extends $\omega|_{M_{\eps}}$. { Using also Weinstein's Lagrangian tubular neighbourhood theorem, we obtain that the Poisson Extension Problem has a positive answer around germs of Lagrangian submanifolds:}

\begin{corollary}
Let $(M^{2n},\omega)$ be a symplectic manifold without boundary and $L\subset M $ be a Lagrangian submanifold. For any neighbourhood $U$ of $L$, there is a Poisson structure $\tilde{\pi}$ on $M$ which has the same germ around $L$ as $\omega$ and is supported in $U$.
\end{corollary}


If $M$ is an orientable four-manifold, then it admits a decomposition $M = M_+ \cup M_-$ where both $M_+$ and $M_-$ are symplectic manifolds with contact boundary \cite{MR2253445}. Since the Poisson structures $\tilde{\pi}_{\pm}$ on $M_{\pm}$ obtained by applying the Poisson Extension Theorem vanish flatly at the boundary, they can be patched to a Poisson structure on $M$. We obtain:

\begin{corollary}
Any orientable four-manifold admits a Poisson structure which is symplectic outside a separating smooth submanifold. 
\end{corollary}

Lanius has developed a method for gluing symplectic manifolds with matching contact boundary into scattering Poisson manifolds \cite{MR4229238}. Lanius' gluing and the results of \cite{MR2253445} also imply the corollary above. These structures are different from the ones produced here.

\subsection{Cosymplectic boundary}\label{ex:cosymplectic}

Let $(M^{2n},\omega)$ be a symplectic manifold with cosymplectic boundary, i.e., there exists an outward pointing vector field $X$, defined around $\partial M$, such that $\mc{L}_X \omega = 0$. Then $\partial M$ is pseudoconvex of regular Pfaffian type of class $1$. Indeed, $\gamma = i_X \omega|_{\partial M}$ is a nowhere vanishing closed 1-form, hence defines a regular Pfaffian distribution of class $1$, $\ka=\ker \gamma$ on $\partial M$. Since its kernel is $\h=\ka$, $\sigma:= \omega|_{\partial M}$ is non-degenerate on the leaves of $\ka$. Hence the first two conditions in Lemma \ref{lem:alternative} hold. Further, since $\mc{L}_X \omega = 0$ the third condition holds for any almost complex structure taming $\omega|_{\ka}$. Hence any boundary component of cosymplectic  type is also of regular Pfaffian type of class $1$. 

{
Therefore, the Poisson Extension Theorem implies also an alternative approach to the second part of Corollary \ref{corollary:contact:and:cosymplectic}. This allows us to conclude that the Poisson Extension Problem holds also along the singular locus of log-symplectic structure:

\begin{corollary}
Let $(M,\pi)$ be a log-symplectic structure with compact singular locus $F$. For any neighbourhood $U$ of $F$, there is a Poisson structure $\tilde{\pi}$ on $M$ which has the same germ around $F$ as $\pi$ and is supported in $U$.
\end{corollary}

\begin{proof}
By the local structure around the singular locus of a log-symplectic structure \cite{MR3250302}, it follows that the boundary of a tubular neighbourhood $T\subset U$ of $F$ is of cosymplectic type. The result follows from the  Poisson Extension Theorem.  
\end{proof}

}

\bibliographystyle{hyperamsplain-nodash}
\bibliography{references}

{{
  \bigskip
  \footnotesize

  G. R.~Cavalcanti, \textsc{Department of Mathematics, Universiteit Utrecht.}\par\nopagebreak
  \textit{E-mail}: \texttt{gil.cavalcanti@gmail.com.}

  \medskip

  I.~M\u{a}rcu\cb{t}, \textsc{Mathematics Department, Radboud University Nijmegen.}\par\nopagebreak
   \textit{E-mail}: \texttt{i.marcut@math.ru.nl}
 }}

\end{document}